\documentclass[12pt]{article}
\usepackage{amssymb, amsmath, amsfonts,dsfont, mathrsfs,euscript,eufrak,colonequals,indentfirst}
\usepackage{textcomp}
\usepackage{latexsym}
\usepackage{ulem}

\newtheorem{Theorem}{Theorem}

\newtheorem{Proposition}{Proposition}

\newcommand{\N}{\mathbb N}
\newcommand{\R}{\mathbb R}
\newcommand{\Z}{\mathbb Z}

\newcommand{\Expec}{ \mathds{E}\,}

\newcommand{\cov}{\mathrm{cov}}

\newcommand{\dd}{\,d}
\newcommand{\id}{\mathds{1}}

\newcommand{\AppClass}{\mathcal{A}}
\newcommand{\CovFunc}{\mathcal{K}}

\newcommand{\e}{\varepsilon}

\newcommand{\tlambda}{\bar{\lambda}}

\newcommand{\tLambda}{\bar{\Lambda}}

\newcommand{\defeq}{\colonequals}
\newcommand{\eqone}{\underset{1}{=}}

\newcommand{\cY}{Z}
\newcommand{\intY}{J}
\newcommand{\appcY}{\widetilde{\cY}}
\newcommand{\appintY}{\widetilde{\intY}}
\newcommand{\SetEP}{\mathcal{S}}
\newcommand{\RF}{V}

\begin{document}





\title{Approximation complexity of homogeneous sums of random processes}
\author{A. A. Khartov, M. Zani}
\maketitle

\begin{abstract}
We study approximation properties of additive random fields $Y_d$, $d\in\N$, which are sums of zero-mean random processes with the same continuous covariance functions. The average case approximation complexity $n^{Y_d}(\e)$ is defined as the  minimal number of evaluations of arbitrary linear functionals needed to approximate $Y_d$, with relative $2$-average error not exceeding a given threshold $\e\in(0,1)$. We investigate the growth of $n^{Y_d}(\e)$ for arbitrary fixed $\e\in(0,1)$ and $d\to\infty$. The results are applied to sums of standard Wiener processes.

\end{abstract}

	

\section{Introduction and problem setting}
Suppose that  we have a sequence of  random processes $\{X_j(t), t\in[0,1], j\in\N\}$ 
defined on some probability space, where $\mathbb N$ denotes the set of positive integers. We assume that all $X_j$ have zero mean and the same continuous covariance function $\CovFunc(t,s)$, $t,s\in[0,1]$. We assume that $\{X_j$, $j\in\N\}$ are uncorrelated. For every $d\in\N$ we define the random field
\begin{eqnarray}\label{def_Yd}
	Y_d(t)\defeq\sum_{j=1}^d  X_j(t_j),\quad t=(t_1,\ldots, t_d)\in[0,1]^d\,.
\end{eqnarray}
This random field has zero mean and covariance function
\begin{eqnarray*}
\CovFunc^{Y_d}(t,s)=\sum_{j=1}^d \CovFunc(t_j,s_j),
\end{eqnarray*}
where $t=(t_1,\ldots, t_d)$ and $s=(s_1,\ldots, s_d)$ are from $[0,1]^d$. Such random fields belong to a wide class of so-called  \textit{additive random fields} (see \cite{ChenLi} and \cite{KarNazNik}).

Every  $Y_d$ is considered as a random element of the space $L_2([0,1]^d)$ endowed with the scalar product $\langle\,\cdot\,,\,\cdot\,\rangle_{2,d}$ and  norm $\|\cdot\|_{2,d}$. We will investigate \textit{the average case approximation complexity} (\textit{approximation complexity} for short) of $Y_d$, $d\in\N$. Recall that the approximation complexity for any random field $\RF_d$ with sample paths from $L_2([0,1]^d)$ is defined by the following formula:
\begin{eqnarray}\label{def_nRFde}
n^{\RF_d}(\e)\colonequals\min\bigl\{n\in\N:\, e^{\RF_d}(n)\leqslant \e\, e^{\RF_d}(0)  \bigr \},
\end{eqnarray}
where $\e\in(0,1)$ is a given error threshold, $e^{\RF_d}(0)\colonequals  \bigl(\Expec\|\RF_d\|_{2,d}^2\bigr)^{1/2}$ is  the total ``size'' of $\RF_d$  ($\Expec$ denotes the expectation), and
\begin{eqnarray*}
e^{\RF_d}(n)\colonequals\inf\Bigl\{ \bigl(\Expec\bigl\|\RF_d - \widetilde \RF_d^{( n)}\bigr\|_{2,d}^2\bigr)^{1/2} :  \widetilde \RF_d^{(n)}\in \AppClass_n^{\RF_d}\Bigr\}
\end{eqnarray*}
is the smallest 2-average  error of approximation  by $n$-rank fields from the class  
\begin{eqnarray*}
	\AppClass_n^{\RF_d}\colonequals \Bigl\{\sum_{m=1}^{n} \langle\RF_d,\psi_m\rangle_{2,d}\,\psi_m :  \psi_m \in L_2([0,1]^d)\Bigr \}.
\end{eqnarray*}
It is well known that $n^{\RF_d}(\e)$ is fully determined by spectral characteristics of the covariance operator $K^{\RF_d}$ of the random field $\RF_d$. Namely, let $(\lambda^{\RF_d}_k)_{k\in\N}$ and $(\psi^{\RF_d}_k)_{k\in\N}$ denote the non-increasing sequence of eigenvalues and the corresponding  sequence of eigenvectors of $K^{\RF_d}$, respectively. Then the following random field
\begin{eqnarray*}
	\RF^{(n)}_d(t)\colonequals\sum_{k=1}^n \langle \RF_d,\psi^{\RF_d}_k\rangle_{2,d}\, \psi^{\RF_d}_k(t),\qquad t\in[0,1]^d,
\end{eqnarray*}
is optimal $n$-rank approximation of $\RF_d$ (see \cite{Rit} and \cite{WasWoz}), i.e.
\begin{eqnarray*}
	e^{\RF_d}(n)=\sqrt{\Expec 
		\bigl\|\RF_d-\RF^{(n)}_d\bigr\|_{2,d}^2},\quad n\in\N.
\end{eqnarray*}
Hence formula \eqref{def_nRFde} is reduced to
\begin{eqnarray*}
	n^{\RF_d}(\e)=\min\Bigl\{n\in\N:\, \Expec 
	\bigl\|\RF_d-\RF^{(n)}_d\bigr\|_{2,d}^2\leqslant\e^2\,\Expec \|\RF_d\|_{2,d}^2 \Bigr\},\quad  \e\in(0,1),
\end{eqnarray*}
and the approximation complexity  $n^{\RF_d}(\e)$ can be described in terms of eigenvalues of $K^{\RF_d}$  (see \cite{NovWoz1}):
\begin{eqnarray*}
n^{\RF_d}(\e)=\min\Bigl\{n\in\N:\, \sum_{k=n+1}^\infty \lambda^{\RF_d}_k\leqslant\e^2\,\Lambda^{\RF_d} \Bigr\},
\end{eqnarray*}
where $\Lambda^{\RF_d}=\Expec \|\RF_d\|_{2,d}^2=\sum_{k\in\N}\lambda^{V_d}_k$ is the trace of $K^{V_d}$.

In the paper we investigate the growth of approximation complexity for additive random fields $Y_d$, $d\in\N$, for arbirarily small fixed error threshold $\e\in(0,1)$ and $d\to\infty$. We assume that all spectral characteristics of marginal covariance $\CovFunc(t,s)$, $t,s\in[0,1]$, or of its given transformations are known. In particular, eigenvalues $\lambda_k$ and corresponding eigenvectors $\psi_k$, are assumed to be known for the covariance operator corresponding to $\CovFunc(t,s)$, $t,s\in[0,1]$. We always set that $(\lambda_k)_{k\in\N}$ is non-increasing. The main difficulty of the problem is that  $\lambda^{Y_d}_k$,  $k\in\N$, are generally unknown or do not easily depend on $\lambda_k$, for large $d\in\N$.

For more general additive random fields the similar problems have been considered  in various settings by Lisfhits and Zani \cite{LifZani1} and \cite{LifZani2}, by Wasilkowski and Wo\'zniakowski \cite{WasWoz2} and \cite{Hick} (with Hickernell), by Khartov and Zani \cite{KhartZani}. In particular, in \cite{KhartZani} an exact formula for $n^{Y_d}(\e)$ was obtained. It should be noted, however, that all these works deal only with the case when covariance operators of marginal processes (the summands of the additive random field) have identical $1$ as an  eigenvector. We are not aware any general results when $1$ is not an eigenvector of $K^{X_j}$. Nevertheless, there exist important processes, which do not satisfy this assumption. The most famous example is the standard Wiener process. The present work is motivated by these facts and is devoted to the  approximation of $Y_d$ without this previous assumption.

The paper is organized as follows. In Section 2 we describe in detail a special decomposition of the random fields $Y_d$, $d\in\N$. This decomposition has an independent interest, because it has a lot of nice properties, which, in addition, will be useful for estimation of $n^{Y_d}(\e)$. It also seems that this decomposition or its modifications can help to investigate the approximation complexity of more general additive random fields.  In Section 3 we obtain the estimates and the exact asymptotics for $n^{Y_d}(\e)$ as $d\to\infty$, and we also describe the case, when $n^{Y_d}(\e)$ is bounded as a function of $d$.

Throughout the paper we will use the following unified notation for the covariance characteristics of random processes and fields. Let $\RF(t)$, $t\in[0,1]^n$, be a given zero-mean random process or field with sample paths from $L_2([0,1]^n)$. We will denote by $K^\RF$ and $\CovFunc^\RF$ the covariance operator and the covariance function of $\RF(t)$, $t\in[0,1]^n$, respectively. Let $(\lambda^{\RF}_k)_{k\in\N}$ denote the sequence of eigenvalues and $(\psi^{\RF}_k)_{k\in\N}$ the corresponding  sequence of eigenvectors of $K^{Z}$. The sequence $(\lambda^{\RF}_k)_{k\in\N}$ is assumed to be ranked in non-increasing order. Then $K^{\RF} \psi^{\RF}_k(t)=\lambda^{\RF}_k\psi^{\RF}_k(t)$, $k\in\N$, $t\in[0,1]^n$. Let $\Lambda^{\RF}$ denote the trace of $K^{\RF}$, i.e. $\Lambda^{\RF}\defeq \sum_{k=1}^\infty \lambda^{\RF}_k$. If $K^\RF$ is of rank $p\in\N$, then we formally set $\lambda_k^{\RF}\colonequals0$ for $k>p$. We denote by $\SetEP(\RF)$ the set of eigenpairs with non-zero eigenvalues:
\[\SetEP(\RF)\defeq\{(\lambda^{\RF}_k,\psi^{\RF}_k):  \lambda^{\RF}_k\ne 0,\, k\in\N\}\,.\] 

Furthermore, we denote by $\R$ the set of real numbers. If $v\in\R^n$, then $v_l$ always denotes $l$-th coordinate of $v$, $l=1,\ldots, n$. For any function $f$ we will denote by $f^{-1}$ the generalized inverse function $f^{-1}(y)\colonequals\inf\bigl\{x\in\R: f(x)\geqslant y\bigr\}$, where $y$ is from the range of $f$. By \textit{distribution function} $F$ we mean a non-decreasing function $F$ on $\R$ that is right-continuous on $\R$, $\lim_{x\to-\infty} F(x)=0$, and $\lim_{x\to\infty} F(x)=1$. For real sequences $a_n$ and $b_n$, $n\in\N$, the relation $a_n\sim b_n$ means that $a_n/b_n\to 1$, $n\to\infty$. The quantity $\id(A)$ equals one for the true relation $A$ and zero for the false one. The number of elements of a finite set $B$ is denoted by $\#(B)$, where $\#(\varnothing)\defeq0$. The function $x\mapsto\lceil x\rceil$ is a ceiling function, i.e.  $\lceil x\rceil=k\in\Z$  whenever $k-1< x\leqslant k$. For any numbers $x$ and $y$ the notation $x\eqone y$ means that $y\leqslant x\leqslant y+1$.

\section{Decomposition of random fields $Y_d$, $d\in\N$}
Let us consider the sequence of random fields $\{Y_d,d\in\N\}$, defined by formula \eqref{def_Yd}. Set $I_j\defeq \int_{[0,1]}X_j(s)\dd s$ and let us define by
\begin{eqnarray*}
	\intY_d\defeq \int\limits_{[0,1]^d} Y_d(s)\dd s=\sum_{j=1}^{d} I_j,
\end{eqnarray*}
and
\begin{eqnarray*}
	\cY_d(t)\defeq Y_d(t)-\int\limits_{[0,1]^d} Y_d(s)\dd s= \sum\limits_{j=1}^{d}(X_j(t_j)-I_j),\quad t\in[0,1]^d.
\end{eqnarray*}
Thus we have
\begin{eqnarray}\label{eq_Decom1}
	Y_d(t)=\intY_d+\cY_d(t),\quad t\in[0,1]^d,\quad d\in\N.
\end{eqnarray}
The summands of decomposition \eqref{eq_Decom1} are orthogonal in $L_2([0,1]^d)$. Every random variable $J_d$ can be considered as one-rank random field whose covariance operator has eigenvector $1$ with associated eigenvalue $\Expec\! J_d^2$. Every random field $Z_d$ is additive and the covariance operators of its marginal processes have identical $1$ as an eigenvector with zero eigenvalue. It is natural to expect that $n^{Y_d}(\e)$ has a similar growth as $n^{Z_d}(\e)$. If such assertion is true, we could apply  results from \cite{KhartZani} to $n^{Z_d}(\e)$. In general, however, $\intY_d$ and $\cY_d$ are correlated, and it is difficult to obtain a connection between $n^{Y_d}(\e)$ and  $n^{Z_d}(\e)$. 
 
We now propose another decomposition for $Y_d$, $d\in\N$, where the parts are orthogonal and uncorrelated and, as we will see below, they are respectively close to $\intY_d$ and $\cY_d$ with small relative errors for large $d\in\N$. We define
\begin{eqnarray*}
	\appintY_d(t)\defeq \sum_{j=1}^d\dfrac{\sum_{l=1}^{d} X_j(t_l)}{d},\qquad \appcY_d(t)\defeq \sum_{j=1}^d \biggl(X_j(t_j)-\dfrac{\sum_{l=1}^{d} X_j(t_l)}{d}\biggr),\quad t\in[0,1]^d.
\end{eqnarray*}
Therefore we have
\begin{eqnarray}\label{eq_decom2}
	Y_d(t)=\appintY_d(t)+\appcY_d(t), \quad t\in[0,1]^d,\quad d\in\N.
\end{eqnarray}

\begin{Proposition}\label{pr_covorth_appintcYd}
For every $d\in\N$ decomposition \eqref{eq_decom2} has the following properties:
\begin{enumerate}
	\item[$1)$] $\cov (\appintY_d(t),\appcY_d(s))=0$,\quad $t,s\in[0,1]^d$;
	\item[$2)$] $\langle\appintY_d,\appcY_d\rangle_{2,d}=0$.
\end{enumerate}	
\end{Proposition}
\textbf{Proof.}\quad Fix $d\in\N$. We first check $1)$. Choose any $t,s\in[0,1]^d$. We have
\begin{eqnarray*}
	\cov(\appintY_d(t),\appcY_d(s))&=&\sum_{j=1}^d \cov \biggl(\dfrac{\sum_{l=1}^{d} X_j(t_l)}{d}, X_j(s_j) -\dfrac{\sum_{k=1}^{d} X_j(s_k)}{d}\biggr)\\
	&=&\dfrac{1}{d}\sum_{j=1}^d \cov \Bigl(\sum_{l=1}^{d} X_j(t_l), X_j(s_j)\Bigr)- \dfrac{1}{d^2}\sum_{j=1}^d\cov \Bigl(\sum_{l=1}^{d} X_j(t_l),\sum_{k=1}^{d} X_j(s_k)\Bigr)\\
	&=&\dfrac{1}{d}\sum_{j=1}^d \sum_{l=1}^{d}\cov \bigl( X_j(t_l), X_j(s_j)\bigr)- \dfrac{1}{d^2}\sum_{j=1}^d\sum_{l=1}^{d}\sum_{k=1}^{d}\cov \bigl( X_j(t_l), X_j(s_k)\bigr).
\end{eqnarray*}
Since $\cov \bigl( X_j(t_l), X_j(s_k)\bigr)=\CovFunc(t_l,s_k)$, we obtain
\begin{eqnarray*}
	\cov(\appintY_d(t),\appcY_d(s))&=&\dfrac{1}{d}\sum_{j=1}^d \sum_{l=1}^{d}\CovFunc(t_l,s_j)- \dfrac{1}{d^2}\sum_{j=1}^d\sum_{l=1}^{d}\sum_{k=1}^{d}\CovFunc (t_l, s_k)\\
	&=&\dfrac{1}{d}\sum_{j=1}^d \sum_{l=1}^{d}\CovFunc(t_l,s_j)- \dfrac{1}{d}\sum_{l=1}^{d}\sum_{k=1}^{d}\CovFunc (t_l, s_k)\\
	&=&0.
\end{eqnarray*}

We now check $2)$. Observe that
\begin{eqnarray*}
	\langle \appintY_d,\appcY_d\rangle_{2,d}=\langle \appintY_d,Y_d-\appintY_d\rangle_{2,d}=\langle \appintY_d,Y_d\rangle_{2,d}-\|\appintY_d\|_{2,d}^2.
\end{eqnarray*}
For $\langle \appintY_d,Y_d\rangle_{2,d}$ we have
\begin{eqnarray*}
	\langle\appintY_d,Y_d\rangle_{2,d}&=&\dfrac{1}{d}\int\limits_{[0,1]^d}\Bigl(\sum_{j=1}^d \sum_{l=1}^{d}X_j(t_l) \cdot \sum_{k=1}^{d} X_k(t_k)\Bigr) \dd t\\
	&=& \dfrac{1}{d}\sum_{j=1}^d\sum_{k=1}^{d}\iint\limits_{[0,1]^2}\Bigl( \sum_{l=1}^{d}X_j(t_l)X_k(t_k)\Bigr) \dd t_l \dd t_k\\
	&=& \dfrac{1}{d}\sum_{j=1}^d\sum_{k=1}^{d}\biggl(\int\limits_{[0,1]} X_j(t_k)X_k(t_k) \dd t_k +\sum_{\substack{l=1,\ldots, d\\ l\ne k}}\,\iint\limits_{[0,1]^2} X_j(t_l)X_k(t_k) \dd t_l \dd t_k \biggr)\\
	&=& \dfrac{1}{d}\sum_{j=1}^d\sum_{k=1}^{d}\biggl(\int\limits_{[0,1]} X_j(t)X_k(t) \dd t +(d-1)I_j I_k \biggr).
\end{eqnarray*}
For $\|\appintY_d\|_{2,d}^2$ we find
\begin{eqnarray*}
	\|\appintY_d\|_{2,d}^2&=&\dfrac{1}{d^2} \int\limits_{[0,1]^d} \Bigl(\sum_{j=1}^{d}\sum_{l=1}^{d}X_j(t_l)\Bigr) \Bigl(\sum_{k=1}^{d}\sum_{r=1}^{d}X_k(t_r)\Bigr)\dd t\\
	&=&\dfrac{1}{d^2}\sum_{j=1}^{d} \sum_{k=1}^{d}\int\limits_{[0,1]^d} \Bigl(\sum_{l=1}^{d}X_j(t_l)\Bigr) \Bigl(\sum_{r=1}^{d}X_k(t_r)\Bigr)\dd t\\
	&=&\dfrac{1}{d^2}\sum_{j=1}^{d} \sum_{k=1}^{d} \Bigl(\sum_{l=1}^{d}\int\limits_{[0,1]}X_j(t_l)X_k(t_l)\dd t_l+\sum_{\substack{l,r=1,\ldots, d\\ l\ne r}}\,\iint\limits_{[0,1]^2}X_j(t_l)X_k(t_r)\dd t_l \dd t_r\Bigr)\\
	&=&\dfrac{1}{d^2}\sum_{j=1}^{d} \sum_{k=1}^{d} \Bigl(d \int\limits_{[0,1]}X_j(t)X_k(t)\dd t+(d^2-d)I_j I_k\Bigr)\\
	&=&\dfrac{1}{d}\sum_{j=1}^{d} \sum_{k=1}^{d} \Bigl( \int\limits_{[0,1]}X_j(t)X_k(t)\dd t+(d-1)I_j I_k\Bigr).
\end{eqnarray*}
We see that $\langle \appintY_d,Y_d\rangle_{2,d}=\|\appintY_d\|_{2,d}^2$, i.e. $\langle \appintY_d,\appcY_d\rangle_{2,d}=0$.\quad$\Box$\\

In Proposition \ref{pr_covorth_appintcYd} the random field $\appintY_d(t)$, $t\in[0,1]^d$, can be replaced by $\intY_d$.
\begin{Proposition}\label{pr_covorth_intappcYd}
For every $d\in\N$ the following properties of $\appcY_d$ hold:
\begin{enumerate}
	\item[$1)$] $\cov (\intY_d,\appcY_d(s))=0$,\quad $t,s\in[0,1]^d$;
	\item[$2)$] $\langle\intY_d,\appcY_d\rangle_{2,d}=\langle 1,\appcY_d\rangle_{2,d}=0$.
\end{enumerate}	
\end{Proposition}
\textbf{Proof.} Fix $d\in\N$. We first check $1)$. Let us choose any $t\in[0,1]^d$. Note that for any $j,l\in\N$
\begin{eqnarray*}
	\cov(I_j, X_l(t))=\Expec \biggl(\,\int\limits_{[0,1]}X_j(s)\dd s \cdot X_l(t)  \biggr)= \int\limits_{[0,1]}\Expec \bigl(X_j(s)X_l(t)\bigr)\dd s.
\end{eqnarray*}
This equals $\int_{[0,1]}\CovFunc(t,s)\dd s$ if $l=j$, and zero if not. Therefore
\begin{eqnarray*}
	\cov(\intY_d,\appcY_d(t))
	&=& \cov\Biggl(\sum_{j=1}^dI_j,\sum_{j=1}^d\biggl( X_j(t_j) -\dfrac{\sum_{k=1}^{d} X_j(t_k)}{d}\biggr)\Biggr)\\
	&=&\sum_{j=1}^d \cov \biggl(I_j, X_j(t_j) -\dfrac{\sum_{k=1}^{d} X_j(t_k)}{d}\biggr)\\
	&=&\sum_{j=1}^d \cov \bigl(I_j, X_j(t_j)\bigr)- \dfrac{1}{d}\sum_{k=1}^{d}\sum_{j=1}^d\cov \bigl(I_j, X_j(t_k)\bigr)\\
	&=&\sum_{j=1}^d \int\limits_{[0,1]}\CovFunc(t_j,s)\dd s- \dfrac{1}{d}\sum_{k=1}^{d}\sum_{j=1}^d\int\limits_{[0,1]}\CovFunc(t_k,s)\dd s.
\end{eqnarray*}
In the last double sum the summands do not depend on $j$. Consequently, we have 
\begin{eqnarray*}
	\cov(\intY_d,\appcY_d(t))=\sum_{j=1}^d \int\limits_{[0,1]}\CovFunc(t_j,s)\dd s- \sum_{k=1}^{d}\int\limits_{[0,1]}\CovFunc(t_k,s)\dd s=0.
\end{eqnarray*}

We now check $2)$:
\begin{eqnarray*}
	\langle 1,\appcY_d\rangle_{2,d}&=& \int\limits_{[0,1]^d} \sum_{j=1}^d \biggl(X_j(s_j)-\dfrac{\sum_{l=1}^d X_j(s_l)}{d} \biggr) \dd s\\
	&=&\sum_{j=1}^d \int\limits_{[0,1]} X_j(s_j) \dd s_j-\dfrac{1}{d}\sum_{j=1}^d \sum_{l=1}^d \int\limits_{[0,1]}X_j(s_l)\dd s_l\\
	&=&\sum_{j=1}^d I_j-\dfrac{1}{d}\sum_{j=1}^d \sum_{l=1}^d I_j\\
	&=&0.
\end{eqnarray*}
Obviously, $\langle \intY_d,\appcY_d\rangle_{2,d}=\intY_d\cdot\langle 1,\appcY_d\rangle_{2,d}=0$.\quad $\Box$\\

In the following proposition we establish an important property concerning eigenpairs of covariance operators in the decomposition \eqref{eq_decom2}.

\begin{Proposition}\label{pr_SetEP}
For every $d\in\N$ the following properties hold:
\begin{enumerate}
\item[$1)$] $\SetEP(\appintY_d)\cap\SetEP(\appcY_d)=\varnothing$;
\item[$2)$] $\SetEP(Y_d)=\SetEP(\appintY_d)\cup\SetEP(\appcY_d)$.
\end{enumerate}
\end{Proposition}
\textbf{Proof.}\quad  Fix $d\in\N$. According to property $1)$ of Proposition \ref{pr_covorth_appintcYd} we have the decomposition
\begin{eqnarray}\label{eq_CovFuncYdappintcYd}
	\CovFunc^{Y_d}(t,s)=\CovFunc^{\appintY_d}(t,s)+\CovFunc^{\appcY_d}(t,s),\quad t,s\in[0,1]^d.
\end{eqnarray}
Let us find the first summand:
\begin{eqnarray*}
	\CovFunc^{\appintY_d}(t,s)&=&\Expec\biggl(\sum_{j=1}^{d} \dfrac{\sum_{l=1}^{d}X_j(t_l)}{d}\cdot \sum_{j=1}^{d} \dfrac{\sum_{r=1}^{d}X_j(s_r)}{d}   \biggr)\\
	&=&\dfrac{1}{d^2}\sum_{j=1}^{d}\Expec \Bigl(\sum_{l=1}^{d}X_j(t_l)\cdot \sum_{r=1}^{d}X_j(s_r)\Bigr)\\
	&=&\dfrac{1}{d^2}\sum_{j=1}^{d}\sum_{l=1}^{d}\sum_{r=1}^{d}\Expec \bigl(X_j(t_l)\cdot X_j(s_r)\bigr)\\
	&=&\dfrac{1}{d^2}\sum_{j=1}^{d}\sum_{l=1}^{d}\sum_{r=1}^{d}\CovFunc(t_l,s_r)\\
	&=&\dfrac{1}{d}\sum_{l=1}^{d}\sum_{r=1}^{d}\CovFunc(t_l,s_r).
\end{eqnarray*}
We represent every $\CovFunc(t_l,s_r)$ by absolutely and uniform convergent series $\sum_{k\in\N}\lambda_k \psi_k(t_l)\psi_k(s_r)$:
\begin{eqnarray*}
	\CovFunc^{\appintY_d}(t,s)&=&\dfrac{1}{d}\sum_{l=1}^{d}\sum_{r=1}^{d}\sum_{k\in\N}\lambda_k \psi_k(t_l)\psi_k(s_r)\\
	&=&\dfrac{1}{d}\sum_{k\in\N}\lambda_k\Bigl(\sum_{l=1}^{d}\psi_k(t_l)\Bigr)\Bigl(\sum_{r=1}^{d} \psi_k(s_r)\Bigr).
\end{eqnarray*}
Suppose that $(\lambda^*,\psi^*)\in\SetEP(\appintY_d)$. Then
\begin{eqnarray}\label{eq_psistar}
	\psi^*(t)=\dfrac{1}{\lambda^*} \int\limits_{[0,1]^d} \CovFunc^{\appintY_d}(t,s)\psi^*(s)\dd s=\sum_{k\in\N}a_{d,k}\Bigl(\sum_{l=1}^{d}\psi_k(t_l)\Bigr),\quad t,s\in[0,1]^d,
\end{eqnarray}
where
\begin{eqnarray*}
	a_{d,k}\defeq \dfrac{1}{d}\cdot\dfrac{\lambda_k}{\lambda^*}\int\limits_{[0,1]^d} \Bigl(\sum_{r=1}^{d} \psi_k(s_r)\Bigr) \psi^*(s)\dd s,\quad k\in\N.
\end{eqnarray*}
We now show that
\begin{eqnarray}\label{eq_CovFuncYdappintYdsum}
	\int\limits_{[0,1]^d}\CovFunc^{Y_d}(t,s) \Bigl(\sum_{l=1}^{d}\psi_m(s_l)\Bigr)\dd s=\int\limits_{[0,1]^d}\CovFunc^{\appintY_d}(t,s) \Bigl(\sum_{l=1}^{d}\psi_m(s_l)\Bigr)\dd s,\quad t\in[0,1]^d.
\end{eqnarray}
Let us denote by $L_d(t)$ and $R_d(t)$ the left-hand side and the right-hand side of \eqref{eq_CovFuncYdappintYdsum}, respectively. Let us consider $L_d(t)$:
\begin{eqnarray*}
	L_d(t)
	&=&\int\limits_{[0,1]^d}\sum_{j=1}^{d}\CovFunc(t_j,s_j)\Bigl(\sum_{l=1}^{d}\psi_m(s_l)\Bigr)\dd s\\
	&=&\sum_{j=1}^{d}\biggl(\sum_{l=1}^{d}\int\limits_{[0,1]^d}\CovFunc(t_j,s_j)\psi_m(s_l)\dd s\biggr)\\
	&=&\sum_{j=1}^{d}\biggl(\int\limits_{[0,1]}\CovFunc(t_j,s_j)\psi_m(s_j)\dd s_j+\sum_{\substack{l=1,\ldots,d\\l\ne j}}\,\iint\limits_{[0,1]^2}\CovFunc(t_j,s_j)\psi_m(s_l)\dd s\dd s_l\biggr)\\
	&=&\sum_{j=1}^{d}\biggl(\int\limits_{[0,1]}\CovFunc(t_j,s)\psi_m(s)\dd s+(d-1)\int\limits_{[0,1]}\CovFunc(t_j,s)\dd s\cdot\int\limits_{[0,1]}\psi_m(s) \dd s\biggr).
\end{eqnarray*}
We next consider $R_d(t)$:
\begin{eqnarray*}
	R_d(t)
	&=& \dfrac{1}{d}\int\limits_{[0,1]^d}\sum_{l=1}^{d}\sum_{r=1}^{d}\CovFunc(t_l,s_r)\Bigl(\sum_{j=1}^{d}\psi_m(s_j)\Bigr)\dd s\\
	&=&\dfrac{1}{d}\sum_{l=1}^{d}\biggl(\sum_{r=1}^{d}\sum_{j=1}^{d}\int\limits_{[0,1]^d}\CovFunc(t_l,s_r)\psi_m(s_j)\dd s\biggr)\\
	&=&\dfrac{1}{d}\sum_{l=1}^{d}\biggl(\sum_{j=1}^{d}\int\limits_{[0,1]}\CovFunc(t_l,s_j)\psi_m(s_j)\dd s_j+\sum_{\substack{j,r=1,\ldots,d\\j\ne r}}\,\iint\limits_{[0,1]^2}\CovFunc(t_l,s_r)\psi_m(s_j)\dd s_r\dd s_j\biggr)\\
	&=&\dfrac{1}{d}\sum_{l=1}^{d}\biggl(d\int\limits_{[0,1]}\CovFunc(t_l,s)\psi_m(s)\dd s+(d^2-d)\int\limits_{[0,1]}\CovFunc(t_l,s)\dd s\cdot\int\limits_{[0,1]}\psi_m(s) \dd s\biggr)\\
	&=&\sum_{l=1}^{d}\biggl(\int\limits_{[0,1]}\CovFunc(t_l,s)\psi_m(s)\dd s+(d-1)\int\limits_{[0,1]}\CovFunc(t_l,s)\dd s\cdot\int\limits_{[0,1]}\psi_m(s) \dd s\biggr).
\end{eqnarray*}
Thus we see that $L_d(t)=R_d(t)$, i.e. \eqref{eq_CovFuncYdappintYdsum} is valid. According to \eqref{eq_psistar}, we conclude that
\begin{eqnarray*}
\int\limits_{[0,1]^d} \CovFunc^{Y_d}(t,s)\psi^*(s)\dd s=\int\limits_{[0,1]^d} \CovFunc^{\appintY_d}(t,s)\psi^*(s)\dd s=\lambda^*\psi^*(t),\quad t\in[0,1]^d.
\end{eqnarray*}
Therefore $(\lambda^*,\psi^*)\in\SetEP(Y_d)$, i.e. $\SetEP(\appintY_d)\subset\SetEP(Y_d)$. Due to \eqref{eq_CovFuncYdappintcYd} we have
\begin{eqnarray*}
	\int\limits_{[0,1]^d} \CovFunc^{\appcY_d}(t,s)\psi^*(s)\dd s=\int\limits_{[0,1]^d} \bigl(\CovFunc^{Y_d}(t,s)-\CovFunc^{\appintY_d}(t,s)\bigr)\psi^*(s)\dd s=0,\quad t\in[0,1]^d.
\end{eqnarray*}
This means that $(\lambda^*,\psi^*)\notin\SetEP(\appcY_d)$. Thus $\SetEP(\appintY_d)\cap\SetEP(\appcY_d)=\varnothing$, i.e. $1)$ is proved. 

We next observe that for any $t,s\in[0,1]^d$
\begin{eqnarray*}
	\CovFunc^{\appcY_d}(t,s)&=&\CovFunc^{Y_d}(t,s)-\CovFunc^{\appintY_d}(t,s)\\
	&=&\sum_{(\lambda,\psi)\in\SetEP(Y_d)}\lambda\, \psi(t)\psi(s) - \sum_{(\lambda,\psi)\in\SetEP(\appintY_d)}\lambda\, \psi(t)\psi(s)\\
	&=&\sum_{(\lambda,\psi)\in\SetEP(Y_d)\setminus\SetEP(\appintY_d)}\lambda\, \psi(t)\psi(s).
\end{eqnarray*}
Since in the last sum all $\psi$ are orthonormal, we have $\SetEP(Y_d)\setminus\SetEP(\appintY_d)\subset\SetEP(\appcY_d)$. Next, we observe that
\begin{eqnarray*}
	\sum_{(\lambda,\psi)\in\SetEP(\appcY_d)}\lambda=\int\limits_{[0,1]^d} \CovFunc^{\appcY_d}(s,s)\dd s=\int\limits_{[0,1]^d}\biggl(\sum_{(\lambda,\psi)\in\SetEP(Y_d)\setminus\SetEP(\appintY_d)}\lambda\, \psi(s)^2\biggr)\dd s=\sum_{(\lambda,\psi)\in\SetEP(Y_d)\setminus\SetEP(\appintY_d)}\lambda.
\end{eqnarray*}
Consequently, there are no $(\lambda,\psi)\in\SetEP(\appcY_d)$ such that $(\lambda,\psi)\notin\SetEP(Y_d)\setminus\SetEP(\appintY_d)$. Therefore
$\SetEP(Y_d)\setminus\SetEP(\appintY_d)=\SetEP(\appcY_d)$. This and $\SetEP(\appintY_d)\subset\SetEP(Y_d)$ yield $2)$.\quad $\Box$\\

Let us compute absolute and relative errors of approximation of $\intY_d$ and $\cY_d$ by $\appintY_d$ and $\appcY_d$ respectively. We introduce the constants:
\begin{eqnarray}\label{def_lambda0tLambda}
	\tlambda_0\defeq \iint\limits_{[0,1]^2} \CovFunc(t,s)\dd t\dd s,\qquad \tLambda\defeq\int\limits_{[0,1]} \CovFunc(s,s)\dd s - \iint\limits_{[0,1]^2} \CovFunc(t,s)\dd t\dd s=\Lambda-\tlambda_0.
\end{eqnarray} 
It is easy to check that   $\tlambda_0$ is equal to $\Expec I_j^2$ and $\tLambda$ is the trace of covariance operator of $X_j-I_j$  for every $j\in\N$.

\begin{Proposition}\label{pr_Abserrors}
	\begin{eqnarray*}
	\Expec\|\intY_d - \appintY_d \|_{2,d}^2=\Expec\|\cY_d - \appcY_d \|_{2,d}^2=\tLambda,\quad d\in\N.
	\end{eqnarray*}
\end{Proposition}
\textbf{Proof.}\quad Fix $d\in\N$. The first equality follows from $Y_d(t)= \intY_d+\cY_d(t)=\appintY_d(t)+\appcY_d(t)$, $t\in[0,1]^d$. Observe that
\begin{eqnarray*}
	\Expec\|\intY_d - \appintY_d \|_{2,d}^2&=&\Expec \int\limits_{[0,1]^d}\biggl( \sum_{j=1}^d\Bigl(I_j- \dfrac{1}{d}\sum_{l=1}^{d} X_j(s_l)\Bigr)\biggr)^2 \dd s\\
	&=&  \int\limits_{[0,1]^d}\Expec\biggl( \sum_{j=1}^d\Bigl(I_j- \dfrac{1}{d}\sum_{l=1}^{d} X_j(s_l)\Bigr)\biggr)^2 \dd s\\
	&=&  \int\limits_{[0,1]^d} \sum_{j=1}^d\Expec\Bigl(I_j- \dfrac{1}{d}\sum_{l=1}^{d} X_j(s_l)\Bigr)^2 \dd s\\
	&=&  \sum_{j=1}^d\Expec\int\limits_{[0,1]^d} \Bigl(I_j- \dfrac{1}{d}\sum_{l=1}^{d} X_j(s_l)\Bigr)^2 \dd s.	
\end{eqnarray*}
Next, we have
\begin{eqnarray*}
	\Bigl(I_j- \dfrac{1}{d}\sum_{l=1}^{d} X_j(s_l)\Bigr)^2&=&I^2_j-2 I_j\, \dfrac{1}{d}\sum_{l=1}^{d} X_j(s_l)+\dfrac{1}{d^2}\sum_{l=1}^{d}\sum_{k=1}^{d} X_j(s_l)X_j(s_k)\\
	&=&I^2_j-2 I_j\, \dfrac{1}{d}\sum_{l=1}^{d} X_j(s_l) +\dfrac{1}{d^2}\sum_{\substack{k,l=1,\ldots, d\\ k\ne l}} X_j(s_l)X_j(s_k)+\dfrac{1}{d^2}\sum_{l=1}^{d}  X_j(s_l)^2.
\end{eqnarray*}
Therefore
\begin{eqnarray*}
	\int\limits_{[0,1]^d}\Bigl(I_j- \dfrac{1}{d}\sum_{l=1}^{d} X_j(s_l)\Bigr)^2\dd s&=&I^2_j-2 I_j\, \dfrac{1}{d}\sum_{l=1}^{d} I_j+\dfrac{1}{d^2}\sum_{\substack{k,l=1,\ldots, d\\ k\ne l}} I_j^2+\dfrac{1}{d^2}\sum_{l=1}^{d} \int\limits_{[0,1]} X_j(s_l)^2\dd s_l\\
	&=& I^2_j-2 I_j\, \dfrac{1}{d}\cdot d I_j+\dfrac{1}{d^2}\cdot (d^2-d) I_j^2+\dfrac{1}{d^2}\cdot d \int\limits_{[0,1]} X_j(s)^2\dd s\\
	&=& \dfrac{1}{d}\cdot\int\limits_{[0,1]} X_j(s)^2\dd s-\dfrac{1}{d}\cdot I_j^2.
\end{eqnarray*}
Hence
\begin{eqnarray*}
	\Expec\|\intY_d - \appintY_d \|_{2,d}^2=\sum_{j=1}^d\Expec \int\limits_{[0,1]^d}\Bigl(I_j- \dfrac{1}{d}\sum_{l=1}^{d} X_j(s_l)\Bigr)^2\dd s=\dfrac{1}{d}\sum_{j=1}^d  \biggl(\Expec\int\limits_{[0,1]} X_j(s)^2\dd s -\Expec I_j^2\biggr).
\end{eqnarray*}
Since
\begin{eqnarray*}
	\Expec\int\limits_{[0,1]} X_j(s)^2\dd s=\int\limits_{[0,1]} \CovFunc(s,s)\dd s=\Lambda,\qquad \Expec I_j^2=\tlambda_0,\quad j\in\N,
\end{eqnarray*}
we obtain
\begin{eqnarray*}
	\Expec\|\intY_d - \appintY_d \|_{2,d}^2=\dfrac{1}{d}\sum_{j=1}^d (\Lambda -\tlambda_0)=\Lambda -\tlambda_0=\tLambda.\quad\Box
\end{eqnarray*}

Let us find $\Expec\|\intY_d\|^2_{2,d}$ and $\Expec\|\cY_d\|^2_{2,d}$. For the first we obtain
\begin{eqnarray*}
	\Expec\|J_d\|^2_{2,d}=\Expec J_d^2= \Expec\Bigl(\sum_{j=1}^{d} I_j\Bigr)^2= \sum_{j=1}^{d} \Expec I_j^2= \tlambda_0 d.
\end{eqnarray*}
For the second, due to the orthogonality $\intY_d$ and $\cY_d$, we have :
\begin{eqnarray*}
	\Expec\|\cY_d\|^2_{2,d}=\Expec\|Y_d\|^2_{2,d}-\Expec\|\intY_d\|^2_{2,d}= \tLambda d.
\end{eqnarray*}

From these formulas and Proposition \ref{pr_Abserrors} we directly obtain the following formulas for relative errors of approximation of $\intY_d$ and $\cY_d$ by $\appintY_d$ and $\appcY_d$, respectively.
\begin{Proposition}\label{pr_Relerrors}
	If $\tlambda_0\ne 0$ then
	\begin{eqnarray*}
	\dfrac{\Expec\|\intY_d - \appintY_d \|_{2,d}^2}{\Expec\|\intY_d\|_{2,d}^2}=\dfrac{\tLambda}{\tlambda_0}\cdot\dfrac{1}{d},\quad d\in\N.
	\end{eqnarray*}
	If $\tLambda\ne 0$ then
	\begin{eqnarray*}
	\dfrac{\Expec\|\cY_d - \appcY_d \|_{2,d}^2}{\Expec\|\cY_d\|_{2,d}^2}=\dfrac{1}{d},\quad d\in\N.
	\end{eqnarray*}	
\end{Proposition}

\section{Approximation complexity of $Y_d$, $d\in\N$}
We now turn to the approximation complexity of $Y_d$, $d\in\N$. The constants $\tlambda_0$ and $\tLambda$ from \eqref{def_lambda0tLambda} are assumed to be known. If $\tLambda=0$ then $X_j(t)=I_j$, $t\in[0,1]$, $j\in\N$, and $Y_d(t)=\intY_d$ and we have  $n^{Y_d}(\e)=1$. 
We next consider only the case $\tLambda\ne 0$. Define $\e_0\defeq(\tLambda/\Lambda)^{1/2}\in(0,1]$.
If $\e_0<1$ then for every $\e\in[\e_0,1)$ we have 
\begin{eqnarray*}
	\Expec\|Y_d-\intY_d\|_{2,d}^2=\Expec\|\cY_d\|_{2,d}^2=\tLambda d\leqslant \e^2 \Lambda d=\e^2 \Expec\|Y_d\|_{2,d}^2.
\end{eqnarray*}
This means that $n^{Y_d}(\e)=1$, $\e\in[\e_0,1)$, $d\in\N$. For the case $\e\in(0,\e_0)$ we have the following result.

\begin{Theorem}\label{th_nYdncYd}
For any $\e\in(0,\e_0)$ and $d\in\N$ 
\begin{eqnarray*}
n^{\cY_d}\bigl(\e/\e_0+d^{-1/2}\bigr)-1	\leqslant n^{Y_d}(\e)\leqslant n^{\cY_d}\bigl(\e/\e_0\bigr)+1.
\end{eqnarray*}
\end{Theorem}
\textbf{Proof.}\quad  Fix any $\e\in(0,\e_0)$ and $d\in\N$. We first prove the upper estimate for $n^{Y_d}(\e)$. Let us consider
\begin{eqnarray*}
	e^{Y_d}(n)^2=\inf\Bigl\{\Expec\bigl\|Y_d - \widetilde Y^{( n)}_d\bigr\|_{2,d}^2:  \widetilde Y^{(n)}_d\in \AppClass_n^{Y_d}\Bigr\},\quad n\in\N.
\end{eqnarray*}
The $n$-rank random field
\begin{eqnarray*}
	Y_{d,1}^{(n)}(t)\defeq\langle Y_d,1\rangle_{2,d}\cdot 1+\sum_{k=1}^{n-1} \langle Y_d,\psi^{\cY_d}_k\rangle_{2,d}\, \psi^{\cY_d}_k(t),\quad t\in[0,1]^d,
\end{eqnarray*}
belongs to the class $\AppClass_n^{Y_d}$. Therefore $e^{Y_d}(n)^2\leqslant \Expec \|Y_d-Y_{d,1}^{(n)} \|_{2,d}^2$, $n\in\N$. By definition of $\cY_d$, we have $\langle \cY_d,1\rangle_{2,d}=0$. Hence $\langle 1,\psi^{\cY_d}_k\rangle_{2,d}=0$ for every $k\in\N$ such that $\lambda^{\cY_d}_k\ne 0$. Indeed,
\begin{eqnarray*}
	\langle 1,\psi^{\cY_d}_k\rangle_{2,d}&=&\dfrac{1}{\lambda^{\cY_d}_k}\cdot\langle 1,\lambda^{\cY_d}_k\psi^{\cY_d}_k\rangle_{2,d}\\
	&=&\dfrac{1}{\lambda^{\cY_d}_k}\cdot\langle 1,K^{\cY_d}\psi^{\cY_d}_k\rangle_{2,d}\\
	&=&\dfrac{1}{\lambda^{\cY_d}_k} \int\limits_{[0,1]^d} \biggl(\int\limits_{[0,1]^d} \CovFunc^{\cY_d}(t,s)\psi^{\cY_d}_k(s)\dd s\biggr)\dd t\\
	&=&\dfrac{1}{\lambda^{\cY_d}_k} \int\limits_{[0,1]^d} \biggl(\int\limits_{[0,1]^d} \Expec\bigl(\cY_d(t)\cY_d(s)\bigr)\psi^{\cY_d}_k(s)\dd s\biggr)\dd t\\
	&=&\dfrac{1}{\lambda^{\cY_d}_k}\,  \Expec\int\limits_{[0,1]^d} \biggl(\,\int\limits_{[0,1]^d}\cY_d(t)\dd t\biggr)\cY_d(s)\psi^{\cY_d}_k(s)\dd s\,,
\end{eqnarray*}
and the inner integral is $\langle \cY_d,1\rangle_{2,d}=0$. For the case $\lambda^{\cY_d}_k=0$ we set $\psi^{\cY_d}_k\equiv 0$.
We can see easily now that $\langle Y_d,1\rangle_{2,d}=\intY_d$ and $\langle Y_d,\psi^{\cY_d}_k\rangle_{2,d}=\langle \cY_d,\psi^{\cY_d}_k\rangle_{2,d}$, $k\in\N$. Therefore
\begin{eqnarray*}
	\Expec \bigl\|Y_d-Y_{d,1}^{(n)} \bigr\|_{2,d}^2=\Expec\Bigl\|\cY_d -\sum_{k=1}^{n-1} \langle \cY_d,\psi^{\cY_d}_k\rangle_{2,d}\, \psi^{\cY_d}_k\Bigr\|_{2,d}^2=e^{\cY_d}(n-1)^2,\quad n\in\N.
\end{eqnarray*}
Thus we obtain the inequality $e^{Y_d}(n)\leqslant e^{\cY_d}(n-1)$, $n\in\N$. Note that 
\begin{eqnarray}\label{eq_ecYd_eYd}
	e^{\cY_d}(0)^2=\Expec\|\cY_d\|^2_{2,d}= \tLambda d= \e_0^2 \Lambda d= \e_0^2\, \Expec\|Y_d\|^2_{2,d}=\e_0^2\, e^{Y_d}(0)^2.
\end{eqnarray}
Therefore we have
\begin{eqnarray*}
	n^{Y_d}(\e)&=&\min\{n\in\N: e^{Y_d}(n)\leqslant \e e^{Y_d}(0) \}\\
	&\leqslant& \min\{n\in\N: e^{\cY_d}(n-1)\leqslant \e e^{Y_d}(0) \}\\
	&\leqslant& \min\{n\in\N: e^{\cY_d}(n)\leqslant \e e^{Y_d}(0) \}+1\\
	&=& \min\{n\in\N: e^{\cY_d}(n)\leqslant (\e/\e_0)\, e^{\cY_d}(0) \}+1\\
	&=& n^{Z_d}(\e/\e_0)+1.
\end{eqnarray*}

We now prove the lower estimate for $n^{Y_d}(\e)$. Let us consider the following $n$-rank random field  for arbitrarily fixed $n\in\N$:
\begin{eqnarray*}
	Y^{(n)}_d(t)\colonequals\sum_{k=1}^n \langle Y_d,\psi^{Y_d}_k\rangle_{2,d}\, \psi^{Y_d}_k(t),\quad t\in[0,1]^d.
\end{eqnarray*}

Due to Proposition \ref{pr_SetEP}, we can write
\begin{eqnarray}
	Y^{(n)}_d
	&=&\sum_{\substack{m=1,\ldots,n,\\(\lambda^{Y_d}_m, \psi^{Y_d}_m)\in\SetEP(Y_d) }} \langle Y_d,\psi^{Y_d}_m\rangle_{2,d}\, \psi^{Y_d}_m\nonumber\\
	&=&\sum_{\substack{m=1,\ldots,n,\\(\lambda^{Y_d}_m, \psi^{Y_d}_m)\in\SetEP(\appintY_d)}} \langle Y_d,\psi^{Y_d}_m\rangle_{2,d}\, \psi^{Y_d}_m
	+\sum_{\substack{m=1,\ldots,n,\\(\lambda^{Y_d}_m, \psi^{Y_d}_m)\in\SetEP(\appcY_d) }} \langle Y_d,\psi^{Y_d}_m\rangle_{2,d}\, \psi^{Y_d}_m\nonumber\\
	&=&\sum_{m=1}^{l_n} \langle Y_d,\psi^{\appintY_d}_m\rangle_{2,d}\, \psi^{\appintY_d}_m
	+\sum_{m=1}^{k_n} \langle Y_d,\psi^{\appcY_d}_m\rangle_{2,d}\, \psi^{\appcY_d}_m,\label{eq_lnkn}
\end{eqnarray}
where 
\begin{eqnarray*}
	l_n\defeq\#\bigl(\bigl\{(\lambda^{Y_d}_m, \psi^{Y_d}_m)\in\SetEP(\appintY_d): m=1,\ldots, n \bigr\}\bigr),\\
	k_n\defeq\#\bigl(\bigl\{(\lambda^{Y_d}_m, \psi^{Y_d}_m)\in\SetEP(\appcY_d): m=1,\ldots, n  \bigr\}\bigr).
\end{eqnarray*}
Here $l_n,k_n\in\{0,1,2,\ldots, n\}$ and $l_n+k_n\leqslant n$. If $l_n=0$ or $k_n=0$, then  we set that the corresponding sum in \eqref{eq_lnkn} equals zero. In \eqref{eq_lnkn} the sequences $(\psi^{\appintY_d}_m)_{m\in\N}$ and $(\psi^{\appcY_d}_m)_{m\in\N}$ are assumed to be ranked in the order corresponding to the following of their elements in the sequence $(\psi^{Y_d}_m)_{m\in\N}$. 

By Propositions \ref{pr_covorth_appintcYd} and \ref{pr_SetEP} we have
\begin{eqnarray*}
	\langle \appintY_d,\appcY_d\rangle_{2,d}=\langle \appintY_d,\psi^{\appcY_d}_r\rangle_{2,d}=\langle \psi^{\appintY_d}_m,\appcY_d\rangle_{2,d}=\langle \psi^{\appintY_d}_m,\psi^{\appcY_d}_r\rangle_{2,d}=0,
\end{eqnarray*}
where $m=1,\ldots, l_n$, and $r=1,\ldots,k_n$, $n\in\N$. Therefore
\begin{eqnarray*}
	\bigl\|Y_d-Y^{(n)}_d\bigr\|_{2,d}^2
	&=&\Bigl\| \appintY_d-\sum_{m=1}^{l_n} \langle Y_d,\psi^{\appintY_d}_m\rangle_{2,d}\, \psi^{\appintY_d}_m+\appcY_d-\sum_{m=1}^{k_n} \langle Y_d,\psi^{\appcY_d}_m\rangle_{2,d}\, \psi^{\appcY_d}_m \Bigr\|_{2,d}^2\\
	&=&\Bigl\| \appintY_d-\sum_{m=1}^{l_n} \langle Y_d,\psi^{\appintY_d}_m\rangle_{2,d}\, \psi^{\appintY_d}_m\Bigr\|_{2,d}^2+\Bigl\|\appcY_d-\sum_{m=1}^{k_n} \langle Y_d,\psi^{\appcY_d}_m\rangle_{2,d}\, \psi^{\appcY_d}_m \Bigr\|_{2,d}^2.
\end{eqnarray*}
From this we conclude that
\begin{eqnarray}\label{conc_eYdn}
	e^{Y_d}(n)^2=\Expec 
	\bigl\|Y_d-Y^{(n)}_d\bigr\|_{2,d}^2\geqslant\Expec \bigl\|\appcY_d-\sum_{m=1}^{k_n} \langle Y_d,\psi^{\appcY_d}_m\rangle_{2,d}\,\psi^{\appcY_d}_m \bigr\|_{2,d}^2.
\end{eqnarray}
We next use the triangle inequality and obtain
\begin{eqnarray}\label{conc_sqrt_appcYdineq}
	\sqrt{\Expec \Bigl\|\appcY_d-\sum_{m=1}^{k_n} \langle Y_d,\psi^{\appcY_d}_m\rangle_{2,d}\, \psi^{\appcY_d}_m\Bigr\|_{2,d}^2}&\geqslant& 	\sqrt{\Expec \Bigl\|\cY_d-\sum_{m=1}^{k_n} \langle Y_d,\psi^{\appcY_d}_m\rangle_{2,d}\, \psi^{\appcY_d}_m\Bigr\|_{2,d}^2}\nonumber\\
	&&{}-\sqrt{\Expec \|\cY_d-\appcY_d\|_{2,d}^2}.	
\end{eqnarray}
By Proposition \ref{pr_covorth_intappcYd}, we have $\langle 1,\appcY_d\rangle_{2,d}=0$ and, consequently, $\langle 1,\psi^{\appcY_d}_m\rangle_{2,d}=0$, $m=1,\ldots, k_n$ (see the proof for $\langle 1,\psi^{\cY_d}_m\rangle_{2,d}=0$ above). Hence
\begin{eqnarray*}
	\langle Y_d,\psi^{\appcY_d}_m\rangle_{2,d}=\intY_d\cdot\langle 1 ,\psi^{\appcY_d}_m\rangle_{2,d}+\langle \cY_d,\psi^{\appcY_d}_m\rangle_{2,d}=\langle \cY_d,\psi^{\appcY_d}_m\rangle_{2,d},\quad m=1,\ldots, k_n.
\end{eqnarray*}
Therefore
\begin{eqnarray}\label{conc_cYdineq}
	\Expec \Bigl\|\cY_d-\sum_{m=1}^{k_n} \langle Y_d,\psi^{\appcY_d}_m\rangle_{2,d}\, \psi^{\appcY_d}_m\Bigr\|_{2,d}^2=\Expec \Bigl\|\cY_d-\sum_{m=1}^{k_n} \langle \cY_d,\psi^{\appcY_d}_m\rangle_{2,d}\, \psi^{\appcY_d}_m\Bigr\|_{2,d}^2\geqslant e^{\cY_d}(k_n)^2.
\end{eqnarray}
According to Proposition \ref{pr_Relerrors}, we have
\begin{eqnarray}\label{conc_appcYdcYd}
	\sqrt{\Expec \|\cY_d-\appcY_d\|_{2,d}^2}=\sqrt{\dfrac{1}{d}\, \Expec\|\cY_d\|^2_{2,d}}=d^{-1/2} e^{\cY_d}(0).
\end{eqnarray}
Combining \eqref{conc_eYdn}--\eqref{conc_appcYdcYd}, we obtain the following inequality
\begin{eqnarray}\label{ineq_eYdnecYdkn}
	e^{Y_d}(n)\geqslant e^{\cY_d}(k_n)-d^{-1/2} e^{Z_d}(0),
\end{eqnarray}
which is valid for any $n\in\N$.

We next turn to estimating of $n^{Y_d}(\e)$. On account of \eqref{eq_ecYd_eYd}, we represent this quantity in the following form:
\begin{eqnarray*}
	n^{Y_d}(\e)=\min\bigl\{n\in\N:\, e^{Y_d}(n)\leqslant \e\, e^{Y_d}(0)  \bigr \}=\min\bigl\{n\in\N:\, e^{Y_d}(n)\leqslant (\e/\e_0)\, e^{\cY_d}(0)  \bigr \}.
\end{eqnarray*}
We now apply \eqref{ineq_eYdnecYdkn} to the latter expression:
\begin{eqnarray*}
	n^{Y_d}(\e)&=&\min\bigl\{n\in\N:\,  e^{\cY_d}(k_n)-d^{-1/2} e^{Z_d}(0)\leqslant (\e/\e_0)\, e^{\cY_d}(0)  \bigr \}\\
	&\geqslant&\min\bigl\{n\in\N:\,  e^{\cY_d}(k_n)\leqslant (\e/\e_0+d^{-1/2})\, e^{\cY_d}(0) \bigr \}\\
	&\geqslant&\min\bigl\{k_n\in\N\cup\{0\}:\, n\in\N,\,  e^{\cY_d}(k_n)\leqslant (\e/\e_0+d^{-1/2})\, e^{\cY_d}(0) \bigr \}\\
	&=&\min\bigl\{k\in\N\cup\{0\}:\,   e^{\cY_d}(k)\leqslant (\e/\e_0+d^{-1/2})\, e^{\cY_d}(0)  \bigr \}\\
	&\geqslant&\min\bigl\{k\in\N:\,  e^{\cY_d}(k)\leqslant (\e/\e_0+d^{-1/2})\, e^{\cY_d}(0)\bigr \}-1.
\end{eqnarray*}
The last minimum is  $n^{\cY_d}(\e/\e_0+d^{-1/2})$ by the definition. Thus we obtain the required lower estimate for $n^{Y_d}(\e)$. \quad$\Box$\\

Theorem \ref{th_nYdncYd} establishes a connection between the approximation complexities of $Y_d$ and $\cY_d$. It is easy to check that for every $\cY_d$ covariance operators of its marginal processes $X_j(t_j)-I_j$, $t_j\in[0,1]$, have identical $1$ as an eigenvector with zero eigenvalue. Thus for the approximation complexity of $\cY_d$, $d\in\N$, we can use the formula from \cite{KhartZani} (see the end of Section 3). Let $(\tlambda_k)_{k\in\N}$ be the sequence of eigenvalues of covariance operator of every $X_j-I_j$. As we mentioned above,  $\tLambda=\sum_{k\in\N}\tlambda_k$. We introduce the distribution function
\begin{eqnarray}\label{def_F}
	F(x)\defeq \dfrac{1}{\tLambda}\sum_{k=1}^\infty \tlambda_k \id\bigl(\tlambda_k\geqslant e^{-x}\bigr),\quad x\in\R.
\end{eqnarray}
In our notation the formula from \cite{KhartZani} takes the form 
\begin{eqnarray}\label{eq_ncYd}
	n^{\cY_d}(\e)\eqone\Biggl\lceil d \cdot\tLambda\! \int\limits_{0}^{1-\e^2} \exp\bigl\{F^{-1}(y)\bigr\}\dd y\Biggr\rceil,\quad \e\in(0,1), \quad d\in\N.
\end{eqnarray}
We see that $n^{\cY_d}(\e)\to\infty$ as $d\to\infty$ for every $\e\in(0,1)$. Also here the integral does not depend on $d$ and it is continuous function of $\e$ in $(0,1)$. Therefore from Theorem \ref{th_nYdncYd} and \eqref{eq_ncYd} we conclude the following result.
\begin{Theorem}\label{th_asympnYd}
For any $\e\in(0,\e_0)$ the following asymptotics holds
\begin{eqnarray*}
n^{Y_d}(\e)\sim d \cdot q(\e),\quad d\to\infty,
\end{eqnarray*}
where
\begin{eqnarray*}
	q(\e)\defeq\tLambda\!\!\!\!\! \int\limits_{0}^{1-(\e/\e_0)^2} \exp\bigl\{F^{-1}(y)\bigr\}\dd y, \quad \e\in(0,\e_0).
\end{eqnarray*}
\end{Theorem}

Thus we see that for any fixed $\e\in(0,\e_0)$ the approximation complexity $n^{Y_d}(\e)$ has linear growth on $d$.

Note that the component $q(\e)$ is in fact the value of approximation complexity of $X_1 - I_1$. Indeed, from \eqref{eq_ncYd} and $\cY_1=X_1-I_1$ we have 

\begin{eqnarray}\label{conc_qe}
	n^{X_1-I_1}(\e/\e_0)\eqone\bigl\lceil q(\e)\bigr\rceil,\quad \e\in(0,\e_0).
\end{eqnarray}

\section{Application to sums of Wiener processes}
Suppose that on some probability space we have the sequence of uncorrelated standard Wiener processes $W_j(t)$, $t\in[0,1]$, $j\in\N$. Here $\CovFunc(t,s)=\min\{t,s\}$, $t,s\in[0,1]$. For every $d\in\N$ we define the random field
\begin{eqnarray*}
\mathbb{W}_d(t)\defeq\sum_{j=1}^d  W_j(t_j),\quad t\in[0,1]^d\,.
\end{eqnarray*}
We consider every random field $\mathbb{W}_d(t)$, $t\in[0,1]^d$, as a random element of the space $L_2([0,1]^d)$, and we study the growth of approximation complexity $n^{\mathbb{W}_d}(\e)$ for any fixed $\e\in(0,1)$ and $d\to\infty$. 

Let us recall necessary spectral characteristics. It is well known that covariance operator of standard Wiener process $W(t)$, $t\in[0,1]$, has the following eigenvalues 
\begin{eqnarray*}
	\lambda_k= \dfrac{1}{\pi^2(k-1/2)^2},\quad k\in\N,
\end{eqnarray*}
and the trace $\Lambda=\sum_{k\in\N}\lambda_k =1/2$ (see \cite{GikhmanSkor} p. 189). Let us consider \textit{the centred Wiener process}, which is defined by the formula 
\begin{eqnarray*}
	W^c(t)\defeq W(t)-\int\limits_{[0,1]}W(s)\dd s,\quad t\in[0,1].
\end{eqnarray*}
For this process and its generalizations there exist results concerning small ball deviations and Karhunen-Lo\`eve expansions (see \cite{BegNikOrs}, \cite{Deheu} and also \cite{Donati}). 
The process $W^c$ has the following covarinace function
\begin{eqnarray*}
	\CovFunc^{W^c}(t,s)=\min\{t,s\}+\dfrac{1}{2}\,(t^2+s^2)-t-s +\dfrac{1}{3},\quad t,s\in[0,1].
\end{eqnarray*} 
The eigenpairs for covariance operator of $W^c$ can be found by differentiation of the equality
\begin{eqnarray*}
	\tlambda f(t)=\int\limits_{[0,1]}\CovFunc^{W^c}(t,s) f(s)\dd s,\quad t\in[0,1].
\end{eqnarray*}
Solution of the corresponding boundary-value problem  yields the eigenvalues (see \cite{BegNikOrs} for more details):
\begin{eqnarray*}
	\tlambda_k= \dfrac{1}{\pi^2 k^2},\quad k\in\N.
\end{eqnarray*}
Hence $\tLambda=\sum_{k\in\N}\lambda_k =1/6$ and distribution function \eqref{def_F} has the form
\begin{eqnarray*}
	F(x)=\dfrac{6}{\pi^2}\sum_{k=1}^\infty \dfrac{1}{k^2}\, \id\bigl(\pi^2 k^2\leqslant e^{x}\bigr),\quad x\in\R.
\end{eqnarray*}

We now turn to the approximation complexity $n^{\mathbb{W}_d}(\e)$. We apply  our general results with $Y_d=\mathbb{W}_d$, $d\in\N$, and $X_j=W_j$, $j\in\N$. The centred Wiener processes $W_j^c(t)\defeq W_j(t)-\int_{[0,1]}W_j(s)\dd s$ correspond to $X_j(t)-I_j$, $t\in[0,1]$, $j\in\N$. We have $\e_0=(\tLambda/\Lambda)^{1/2}=3^{-1/2}$. Therefore $n^{\mathbb{W}_d}(\e)=1$ for $\e\in[3^{-1/2},1)$ and $d\in\N$ (see comments before Theorem \ref{th_nYdncYd}). By Theorem \ref{th_asympnYd}, for any $\e\in(0,3^{-1/2})$ we have the asymptotics
\begin{eqnarray*}
	n^{\mathbb{W}_d}(\e)\sim d\cdot q(\e),\quad d\to\infty,
\end{eqnarray*} 
where
\begin{eqnarray*}
	q(\e) \defeq\dfrac{1}{6}\cdot\!\!\!\!\! \int\limits_{0}^{1-3\e^2} \exp\bigl\{F^{-1}(y)\bigr\}\dd y,\quad \e\in (0,3^{-1/2}).
\end{eqnarray*} 

Let us study the behaviour of the component $q(\e)$ for small $\e$. From \eqref{conc_qe} we have
\begin{eqnarray}\label{conc_qeWc}
	n^{W^c_1}(3^{1/2} \e)\eqone \bigl\lceil q(\e)\bigr\rceil, \quad \e\in (0,3^{-1/2}).
\end{eqnarray}
Here
\begin{eqnarray*}
	n^{W^c_1}(3^{1/2} \e)&=&\min\biggl\{n\in\N: \sum_{k=n+1}^\infty \dfrac{1}{\pi^2 k^2} \leqslant 3 \e^2 \cdot \dfrac{1}{6}\biggr\}\\
	&=& \min\biggl\{n\in\N: \sum_{k=n+1}^\infty \dfrac{1}{k^2} \leqslant  \dfrac{\pi^2}{2}\cdot \e^2\biggr\}.
\end{eqnarray*}
Since $\sum_{k=n}^\infty 1/k^2\sim 1/n$, $n\to\infty$, we have $n^{W^c_1}(3^{1/2} \e)\sim (2/\pi^2)\cdot \e^{-2}$, $\e\to 0$. Due to \eqref{conc_qeWc}, we have the asymptotics
\begin{eqnarray*}
	q(\e)\sim \dfrac{2}{\pi^2}\cdot \e^{-2},\quad \e\to 0.
\end{eqnarray*}

\section*{Acknowlegements}
The work of the first named author was supported by DFG--SPbSU grant 6.65.37.2017, the RFBR grant 16-01-00258. Also the work of the first named author was also partially supported by the Government of the Russian Federation (grant 074-U01).

The authors are pleased to thank Institute Denis-Poisson, University of Orleans, where the work was prepared.

\textit{Keywords and phrases}:
additive random fields, average case approximation complexity, asymptotic analysis, sums of Wiener processes.\\

\textbf{A. A. Khartov}\\
\textsc{St. Petersburg State University, 7/9 Universitetskaya nab., 199034 St. Petersburg, Russia.}\\
\textsc{ITMO UNIVERSITY, Kronverksky Pr. 49, 197101 St. Petersburg,  Russia.}\\
\textit{E-mail address}: \texttt{alexeykhartov@gmail.com}\\

\textbf{M. Zani}\\
\textsc{Institut Denis Poisson, Universit\'e d'Orl\'eans,  B\^{a}timent Math\'ematiques, Rue de Chartres, B.P. 6759-45067, Orl\'eans cedex 2, France.}\\
\textit{E-mail address}: \texttt{marguerite.zani@univ-orleans.fr}\\


\begin{thebibliography}{99}
\bibitem{BegNikOrs} L. Beghin, Ya. Nikitin, E. Orsingher, \textit{Exact small ball constants for some Gaussian processes under the $L^2$-norm}, Journal of Math. Sciences, \textbf{128} (2005), 1, 2493--2502.

\bibitem{ChenLi} X. Chen, W.V. Li, \textit{Small deviation estimates for some additive processes}, Proc. Conf. High Dimensional Probab. III, Progress in Probability, \textbf{55} (2003), Birkh\"{a}user, 225--238.

\bibitem{Deheu} P. Deheuvels, \textit{Karhunen-Lo\`eve expansions of mean-centred Wiener processes}, IMS Lecture Notes-Monograph Series, High Dimensional Probab., \textbf{51} (2006), 62--76.

\bibitem{Donati} C. Donati-Martin, M. Yor, \textit{Fubini’s theorem for doubleWiener integrals and the variance of the Brownian path}, Ann. Inst. H. Poincar\'e, \textbf{27} (1991), 181--200.

\bibitem{GikhmanSkor} I. I. Gikhman, A. V. Skorokhod, \textit{Introduction to the Theory of Random Processes}, W.B. Saunders Comp., Philadelphia, 1969.	

\bibitem{Hick} F. J. Hickernell, G. W. Wasilkowski, H. Wo\'zniakowski, \textit{Tractability of linear multivariate problems in the average-case setting}, in: A. Keller, S. Heinrich, H. Niederreiter (Eds.), Monte Carlo and Quasi-Monte Carlo Methods 2006, Springer, Berlin, 2008, pp. 461--493.

\bibitem{KarNazNik} A. Karol, A. Nazarov, Ya. Nikitin, \textit{Small ball probabilities for Gaussian random fields and tensor products of compact operators}, Trans. Amer. Math. Soc., \textbf{360} (2008), no. 3, 1443--1474.

\bibitem{KhartZani} A. A. Khartov, M. Zani,  \textit{Asymptotic analysis of average case approximation complexity of additive random fields}, J. Complexity (Accepted).

\bibitem{LifZani1} M. A. Lifshits, M. Zani, \textit{Approximation complexity of additive random fields}, J. Complexity, \textbf{24} (2008), no. 3, 362--379.

\bibitem{LifZani2} M. A. Lifshits, M. Zani, \textit{Approximation of additive random fields based on standard information: Average case and probabilistic settings}, J. Complexity, \textbf{31} (2015), no. 5, 659--674.

\bibitem{NovWoz1}  E. Novak, H. Wo\'zniakowski, \textit{Tractability of Multivariate Problems. Volume I: Linear Information}, EMS Tracts Math. 6, EMS, Z\"urich, 2008.

\bibitem{Rit} K. Ritter, \textit{Average-case Analysis of Numerical Problems}, Lecture Notes in Math. No. 1733, Springer, Berlin, 2000.

\bibitem{WasWoz} G. W. Wasilkowski, H. Wo\'zniakowski, \textit{Average case optimal algorithms in Hilbert spaces}, J. Approx. Theory, \textbf{47} (1986), 17--25.

\bibitem{WasWoz2} G. W. Wasilkowski, H. Wo\'zniakowski, \textit{Polynomial-time algorithms for multivariate linear problems with finite-order weights: average case setting},  Found. Comput. Math., \textbf{9} (2009), 105--132.
\end{thebibliography}
\end{document}